\documentclass[a4paper]{article}\let\received\date

\received{GF '04, Novi Sad, 27 September 2004}

\title{Generalized functions as sequence spaces with ultranorms}

\author{Maximilian F. \textsc{Hasler}, 
Universit\'e Antilles-Guyane, D.S.I., F--97275 Schoelcher}

\usepackage{amsmath,amsfonts,amssymb,t1enc}
\makeatletter

\catcode176=13 \def^^b0{\degres}
\catcode178=13 \def^^b2{\ensuremath{^2}}
\catcode183=13 \def^^b7{\ensuremath{\cdot}}

\ifx\undefined\doi 
\def\tx#1{\advance\textheight#1 #1=0pt}\tx\topmargin\tx\headsep\tx\headheight
\advance\textwidth2em\advance\hoffset-1em
\def\myfoot{{\sl M.Hasler:~Generalized...ultranorms}\dotfill--\thepage--\dotfill}
\def\@oddfoot{\myfoot(compiled \today)}
\let\mydate\@date
\fi

\ifx\undefined\doi
\newtheorem{theorem}{Theorem}
\numberwithin{theorem}{section}%

\newtheorem{definition}[theorem]{Definition}
\newtheorem{example}[theorem]{Example}
\newtheorem{lemma}[theorem]{Lemma}
\newtheorem{proposition}[theorem]{Proposition}
\newtheorem{remark}[theorem]{Remark}
\newenvironment{proof}[1][Proof]{\nb{#1.}}
{\kern0pt\mbox{}\hfill$\Box$\par}
\fi

\newcommand\A{{\mathcal A}}
\newcommand\B{{\mathcal B}}
\newcommand\C{{\mathbb C}}
\newcommand\CC{{\mathcal C}}
\newcommand\CEP{$\left(\mathcal C,\mathcal E,\mathcal P\right)$}

\newcommand\D{{\mathcal D}}
\newcommand\DD{{\text{\textbf{\textsl D}}}}

\newcommand\E{{\mathcal E}}

\newcommand\F{{\mathcal F}}

\newcommand\G{{\mathcal G}}

\newcommand\ie{\emph{i.e.}}
\def\INT#1d{\int#1\,\mathrm d}
\newcommand\indlim{\operatorname*{ind\,lim}}
\newcommand\impl{\mathop{~\Longrightarrow~}}
\newcommand\J{{\mathcal J}}
\newcommand\K{{\mathcal K}}
\newcommand\KK{{\mathbb K}}
\def\l{\lambda}

\newcommand\lr[3]{\left#1#3\right#2}
\newcommand\Mid{~\Big|~}
\def\myeqn#1{\\[1ex]\centerline{$\begin{array}{r@{}l}\displaystyle#1\end{array}$}\\[1ex]}
\newcommand\N{{\mathbb N}}
\newcommand\nb[1]{\par\noindent\textbf{#1} }
\def\ntoi{\limits_{n\to\infty}}
\def\O{{\mathcal O}}
\newcommand\ola{\overleftarrow}
\newcommand\olra{\!\overleftrightarrow}
\newcommand\ora{\overrightarrow}
\def\P{{\mathcal P}}
\newcommand\p[1]{\left(#1\right)}
\newcommand\pr[1]{_{p,r^{#1}}}

\newcommand\R{{\mathbb R}}
\newcommand\set[1]{\left\{\,#1\,\right\}}
\newcommand\supp{\operatorname*{supp}}
\newcommand\T{{\mathbb T}}
\newcommand\ultra[1]{{|\!|\!|\,#1\,|\!|\!|}} 

\newcommand\vp{\varphi}
\newcommand\veps{\varepsilon}
\newcommand\Z{{\mathbb Z}}

\begin{document}
\maketitle
\begin{abstract}
We review our recent formulation (with A.~Delcroix, S.~Pilipovi\'c and V.~Valmorin)
of Colombeau type algebras as Hausdorff sequence spaces with ultranorms,
defined by sequences of exponential weights.
We extend previous results and give several new perspectives related to echelon type spaces,
possible generalisations, asymptotic algebras, concepts of association, and applications thereof.
\\[1.5ex]
\end{abstract}

\section{Introduction}

Colombeau's New Generalized Functions~\cite{co} are today the most 
widely used 
associative differential algebras containing the $\delta$--distribution.
Their topology is studied since the late 1990s~\cite{infexp}, and
investigation in topological duals of such spaces is now emerging
as important topic of research in this field.

We define such algebras right from the start as spaces with ultranorms~\cite{DHPV1,DHPV3},\linebreak
which is natural and 
especially useful for practical use of the topology,
with no need for valuations.
Our construction allows for algebras containing ultradistributions
and periodic hyperfunctions~\cite{DHPV2}.
Without specializing to a concrete space, 
we deduce general results about completeness, embedding of duals and functoriality,
and generalize known concepts of association, revealing aspects of 
the underlying structure rather hidden in other approaches.
%
Our approach also shows better the close link with the 
classical theory of sequence spaces.

\section{The basic construction}


Consider a sequence $r={(r_n)}_n\in(\R_+)^\N$ decreasing to zero.
For a seminorm $p$ 
on an $\R$-- or $\C$--vector space $E$,
this defines a map $\ultra\cdot_{p,r}:E^\N\to\overline\R_+=[0,\infty]$,\\[1ex]
\centerline{$\displaystyle f={(f_n)}_n~~\longmapsto~~ \ultra{f}=\ultra{f}_{p,r} =
{\limsup\ntoi}\big(p(f_n)\big)^{r_n}~.$}

\begin{lemma}\label{lem1}
~(a)~ If ~$0<\liminf p(f_n)\le\limsup p(f_n)<\infty$, then $\ultra f=1$.\\
(b)~ For all $f,g\in E^\N,\, \l\in\C^*: \ultra{f+g}\le\max(\ultra f,\ultra g)$
and $\ultra{\l f}=\ultra f$.\\
(c)~ If $E$ is a topological algebra,
then $\ultra{f\cdot g}\le\ultra f\cdot\ultra g$.
\end{lemma}

\begin{proof}
As $\lim r_n=0$, we have $\lim k^{r_n}=1$ for any $k>0$, thus (a).
Writing $p(\lambda f_n)\le|\lambda|\,p(f_n)$ and $p(f_n+g_n)\le
2\max(p(f_n),p(g_n))$, we have (b), and using 
$\exists C>0~\forall x,y\in E: p(x\,y)\le C\,p(x)\,p(y)$, we get (c) 
in the same way.
\end{proof}

\begin{definition}
\textbf{The $r$--generalized semi-normed space $(E,p)$}
is the factor space ~$\G_r(E,p) = \F_r(E,p)\,/\,\K_r(E,p)$, where
\\[1ex]\centerline{
$\F_r(E,p) = \big\{ f\in E^\N\mid\ultra{f}_{p,r}<\infty \big\} ~,~~
\K_r(E,p) = \big\{ f\in E^\N\mid\ultra{f}_{p,r} = 0 \big\}~.$}
\end{definition}

\begin{proposition}
The map $\ultra\cdot_{p,r}$ defines a pseudometric $d_{p,r}$ on $\F_r(E,p)$,
making it a topological ring, if $E$ is a topological algebra.
As $\K_r(E,p)$ is the intersection of neighborhoods of zero,
$\G_r(E,p)$ is then the associated Hausdorff topological ring,
on which $d\pr{}$ is well-defined and an ultrametric.
\end{proposition}

\begin{proof}
This is a direct consequence of the Definition and preceding Lemma.
\end{proof}

\begin{example}
For $E=\C$ and $p=|\cdot|$, we obtain
\textbf{the ring of $r$--generalized complex numbers}, $\C_r=\G_r(\C,|·|)$.
%
For $r_n=\frac1{\log n}~(n>1)$, this are 
\textbf{Colombeau's generalized numbers} $\smash{\widetilde\C}$,
since 
\(
\limsup|x_n|^{1/\log n} < \infty
\iff\exists\gamma\in\R: |x_n| = o(n^\gamma),
\) 
and $
\lim|x_n|^{1/\log n} =0\iff\forall\gamma\in\R: |x_n| = o(n^\gamma)$.
\\
The choice $r_n = n^{-1/m}$ (with $m>0$), leads to 
\textbf{ultracomplex numbers} $\overline\C^{p!^m}\!.$
\end{example}


\begin{proposition}\label{top-ring}
The spaces $\G_r(E,p)$ (resp. $\F_r(E,p)$) are topological algebras over
the generalized numbers (resp. over $\F_r(\KK,|\cdot|)$) equipped with
$\ultra·$--topology, but they aren't topological vector spaces
over the field $\KK=\R$ or $\C.$
\end{proposition}

\begin{proof}
This is seen by observing that Lemma~\ref{lem1}-(c) also holds for
$f\in\C^\N$, while  Lemma~\ref{lem1}-(b) implies that $\ultra{\l\,f}$
does not go to zero when $\l\to0$.
\end{proof}

\begin{example}
\label{Ex-Sobolev}\def\W{{W^{s,\infty}(\Omega)}}%
To obtain $r$--\textbf{generalized Sobolev algebras}
$\G_\W=\G_r\big(\W,p_{s,\infty}\big)$, we choose $E=\W$ with norm $p
_{s,\infty}=\smash{\sum\limits_{|\alpha|<s}}{\|\partial^\alpha\cdot\|}_{L^\infty}
$.
This generalizes to any normed algebra.
\end{example}

\begin{theorem}[(equivalent scales)]\label{equiv-scales}
If $r={(r_n)}_n$, $s={(s_n)}_n$ decrease to zero such that
$s=O(r)$, then $\F_r(E,p)\subset\F_s(E,p),\ \K_s(E,p)\subset\K_r(E,p).\,$
In particular, if $ \lim\ntoi\frac{s_n}{r_n}=C\in\R_+^* $,
then $ \ultra f_{p,s} = (\ultra f_{p,r})^C $, and thus
$\F_s(E,p)=\F_r(E,p),\ \K_s(E,p)=\K_r(E,p)$ and $\,\G_s(E,p)=\G_r(E,p).$
\end{theorem}

\begin{proof}
If $s_n=c_n\,r_n$ with $\limsup c_n=C\in\R_+^*$, we have
$\log\ultra f_{p,s} = \limsup (s_n\log p(f_n)) = \limsup (c_n\,r_n\log p(f_n))
\le C\,\limsup (r_n\log p(f_n))= C\log \ultra f_{p,r}$,
where we assumed $\lim\log p(f_n)\ge0$, \ie\ $\ultra f\ge1$.
Otherwise, $\le$ must be replaced by $\ge$,
leading to the inverse inclusion for $\K$.
\end{proof}

\def\avi{_r(\C,\lr||\cdot)}
\begin{remark}[(relation to \textsc{Maddox}' sequence spaces) ]
Our spaces $\K\avi$ and $\F\avi$ are the same as 
$c_0(r)=\bigcap
_{k\in\N}\big\{x\in\C^\N\mid{\lim|x_n|\,k^{1/{r_n}}=0}\big\}$
and\linebreak[3] $\ell_{\infty}(r) =\bigcup
_{k\in\N}\set{x\in\C^\N\mid\smash{\sup|x_n|\,k^{-1/{r_n}}<\infty}}$,
introduced in~\cite{nakano17,simons19} and studied
extensively by Maddox and his students~\cite{mad10,mad14}. 
To see this, observe that %
$\exists k\in\N:\linebreak[0]\sup|x_n|\,k^{-1/r_n}<\infty
	\iff \exists k:\limsup|x_n|^{r_n}\le k \iff \ultra x_r<\infty,$
and ~$
	\forall k:\lim|x_n|\,k^{1/r_n}=0
	\iff \forall \veps>0: |x_n|=o(\veps^{1/r_n})
	\iff \ultra x_r=0
.$\par

These spaces belong to the 
classes of echelon resp.\
co-echelon spaces. 
%
As we always require $\lim r_n=0$, 
both are Montel and Schwartz spaces. 
%
While the cited 
work on sequence spaces is restricted to $(\C,\lr||\cdot)$,
our studies concern 
more general spaces. 
However, most of the spaces considered in the sequel can be written as
intersection and/or union of echelon and co-echelon type spaces.
This also allows the 
generalization of the present construction 
to any abstract topological module $E$, as will be discussed in a forthcoming publication.
\end{remark}

\section{Generalized locally convex spaces}

\begin{definition}
The \textbf{$r$--extension of  a locally convex space} $(E,\P)$ is
the factor space \smash{$ \G_r(E,\P) = \F_r(E,\P)\,/\,\K_r(E,\P)
=\bigcap\limits_{p\in\P}\F_r(E,p)
\,\big/\! 
\bigcap\limits_{p\in\P}\K_r(E,p)
$.}\vskip2ex
\end{definition}

\begin{theorem}
If $(E,\P)$ is a topological algebra,
\ifx\undefined\doi
\ie\ 
\(
  \forall p\in\P~\exists \bar{p}\in\P~\exists C>0\linebreak[3]
  \forall x,y\in E: p(x\,y) \le C\,\bar{p}(x)\,\bar{p}(y) ,
\)
\fi
then $\F_r(E,\P)$ is a subalgebra of $E^\N$,
{$\K_r(E,\P)$}\ is an ideal of $\F_r(E,\P)$,
and $\lr(){d_{p,r}}_{p\in\P}$ is a family of pseudo-distances
on $\G_r(E,\P)$ making it a Hausdoff topological algebra over $\C_r$.
\end{theorem}

\begin{proof}\def\Pr{_r}
Lemma~\ref{lem1}-(b) yields for $f,g\in\F\Pr$, $\lambda\in\C$, and
$p\in\P:\ultra{\lambda f+g}_p\le\max(\ultra f_p,\ultra g_p)$,
thus $\F\Pr$ and $\K\Pr$ are $\C$--linear subspaces.
Continuity of multiplication in $(E,\P)$ gives as in 
\ref{lem1}-(c),
$\forall p\in\P,\ \exists\bar p\in\P:$ $\ultra{f\,g}_p\le\ultra f_{\bar p}\cdot\ultra g_{\bar p}$.
Thus $\F\Pr$ is a $\C$--subalgebra of $E^\N$, and $\K\Pr$ an ideal of $\F\Pr$.
The 
inequalities also imply continuity of addition and multiplication, 
thus $\F\Pr$ is a 
topological $\F_{|\cdot|,r}$--algebra, and
%
$\G_\P$ is again 
the associated Hausdorff space.
\end{proof}

\begin{example}
\label{Ex-simple-Colombeau}
The classical \textbf{simplified Colombeau algebra}~\cite{co}
is obtained for $r_n=\frac{1}{\log n}$ 
and $\P=\set{p_\nu^\mu:f\mapsto\!\!\!\smash{
\sup\limits_{|\alpha|\le\nu,|x|\le\mu}}\! |f^{(\alpha)}(x)|}_{\mu,\nu\in\N }$
on $E=\CC^\infty(\Omega)$.\vskip2ex
\end{example}


\def\projlim{\,\mathop{\!\smash{\rm proj\,lim}}}

As 
a last generalization of the base space, 
consider a 
family of semi-normed algebras
$\left( E_{\nu}^{\mu},p_{\nu\,}^{\mu}\right) _{\mu,\nu\in\N}$ with 
embeddings 
$ \forall\mu,\nu\in\N : E_{\nu}^{\mu+1}\hookrightarrow E_{\nu}^{\mu}$, 
$ E_\nu^\mu \hookrightarrow E_{\nu+1}^\mu$
resp.\ $ E_{\nu+1}^{\mu}\hookrightarrow E_{\nu}^{\mu} $.
Let
$\ora{E} = \projlim\limits_{\,\mu\to\infty} \indlim\limits_{\nu\to\infty} E_\nu^\mu$,
resp.
$\ola{E} = \projlim\limits_{\mu\to\infty} \projlim\limits_{\nu\to\infty}E_\nu^\mu$,
and assume that for all $\mu\in\N$ the inductive limit is regular,
\ie\ a subset is bounded iff it is a bounded subset of
${\left(E_\nu^\mu\right)^\N}$ for some $\nu\in\N$.
%
Now let\\[.5ex]
{$
  \F_r(\ora E) = \set{ f\in\ora E^\N \Mid
  \forall\mu\in\N,\exists\nu\in\N : f\in\lr(){E_\nu^\mu}^\N
  \land \ultra{f}_{p_\nu^\mu,\,r} < \infty } 
~,$}\\
$  \F_{r}(\ola E) =\big\{ f\in\ola{E}^\N \mid
  \forall\mu,\nu\in\N : \ultra{f}_{p_\nu^\mu,\,r}<\infty \big\}
,$ 
and the obvious definition for $\K_r(\olra{E})$, where we
write $\olra\cdot$ for both, $\overrightarrow{\cdot}$ and $\ola\cdot$.
Then again, $\G_r(\olra E)=\F_r(\olra E)\,/\,\K_r(\olra{E})$ is a
topological algebra for the respective 
limit topology.

\begin{example}
\def\rmp{_{r^{m'}}}%
In \cite{DHPV2} we showed how to embed
Gevery class ultradistributions in Colombeau algebras $
  \G^{(p!^m,p!^{m'})}=\G\rmp(\E^{(m)}) ,~
  \G^{\{p!^m,p!^{m'}\}}=\G\rmp(\E^{\{m\}}) ,
$
where $r^m_n=n^{\frac1m}$ and $\E^{(m)},\ \E^{\{m\}}$ are the proj--proj resp.
proj--ind limit type spaces of Beurling resp. Roumieu ultradifferentiable
functions, defined trough spaces on which 
$
    p_\nu^{m,\mu}(f) = \smash{\sup\limits_{|x|\le\mu, \alpha\in\N^s}}
    \dfrac{\nu^{|\alpha|}}{{\alpha!}^m} |f^{(\alpha)}(x)| $ resp. $
    q^{m,\mu}_{\nu}=p^{m,\mu}_{1/\nu} 
$ are finite.
\end{example}


\def\ultrapm#1{\ultra{#1}^{{}_\pm}}
\def\ultrapmk#1{\ultrapm{{(#1)}_k}}

\begin{definition}
\label{sect:periodic}%
Consider the spaces $\O_\l$ 
of holomorphic functions on
$\Omega_\l=\set{z\in\C\mid \frac1\l<|z|<\l}$ with finite norm
$ q^\l=\left\|\cdot\right\|_{L^\infty(\Omega_\l)} $.
Analytic functions on the unit circle are then 
 $\A(\T)=\indlim\limits_{\l\to1}\O_\l$.
Let
$\ora E=\indlim\limits_{\l\to1}(\A_1(\T),q^\lambda)$, where\\\mbox{}\hfil
$\,\displaystyle
   \A_1(\T)=\indlim_{m\to1^-} \indlim_{\nu\to\infty}
   \left\{ f\in \A(\T) \mid
   {\| f^{(\alpha)} \|}_{L^\infty(\T)} \underset{\alpha\to\infty}
   = O({\nu^{\alpha}\,\alpha!^m}) \right\}.
$\\[1ex]
Then, \textbf{$r$--generalized hyperfunctions on $\T$} are defined as
$    \G_{H,r}(\T) = \G_{r}(\ora E) ,
$ 
quotient space of
$
   \F_r(\ora E)=
   \bigcup\limits_{\l>0}\F_r(\A_1(\T),q^\l) 
$ by $
   \K_r(\ora E)=
   \bigcup\limits_{\l>0}\K_r(\A_1(\T),q^\l) .
$
\end{definition}

The same type of ultranorm can be used to characterize generalized
functions $f$ on the unit circle by means of their Fourier coefficients
${(\widehat f_k)}_k\in\C^\Z$, for which we define
$~
  \ultrapm{{(\widehat f_k)}_{k\in\Z}}_r 
  = \max\big\{ \ultra{ {(\widehat f_k)}_{k\in\N} }_{|\cdot|,r}\,,~
  		\ultra{{(\widehat f_{-k})}_{k\in\N}}_{|\cdot|,r} \big\}
.\ $ \par
Fourier coefficients of analytic functions $f\in\A(\T)$,
Schwartz distributions $T\in\E'(\T)$
and hyperfunctions $H\in\B(\T)$ are characterized by
$\ultrapmk{\widehat f_k}_{(\cdot)^{-1}}<1$,
$\ultrapmk{\widehat T_k}_{1/\log}<\infty$, resp.
$\ultrapmk{\widehat H_k}_{(\cdot)^{-1}}\le1$.

\begin{proposition}[(Fourier characterization)]\label{aba}
The same spaces $\F_r(\ora E)$, $\K_r(\ora E)$
are obtained in the previous definition if $q^\l$ is replaced 
by $\widehat{q}^\l:f\mapsto\smash{\displaystyle
\sup\limits_{k\in\Z}}\,{\l^{|k|}|\widehat{f}_k|}$.
\end{proposition}

\begin{proof}
If $f\in\F_r(\ora E)$, $\ultra{f}_{q^{\l},r} <\infty$,
there is $C>0$ such that $q^\l(f_n)^{r_n}<C$ for all $n$.
Cauchy's inequalities in $\Omega_\l$ then give
$|\hat{f}_n(k)|\le q^{\l}(f_n){\l}^{-|k|}$, thus
$|\hat{f}_n(k)|^{r_n}\le C\,\l^{-|k|r_n}$ for all $k\in\Z $,
whence $\ultra{f}_{\hat{q}^\l,r}<\infty$.
Conversely, if $f\in\ora E^\N$, $\ultra{f}_{\hat{q}^{\l},r}<\infty$,
we have $\hat{q}^{\l}(f_n)^{r_n}<C$ for some $C>0$ and all $n$,
\ie\ $|\hat{f}_n(k)|<C^{1/r_n} \l^{-|k|}$ for all $k\in\Z$.
Consequently, there is $M>0$ such that 
$q^\l(f_n) \le M\,C^{1/r_n}$, thus $\ultra{f}_{q^\l,r}<\infty$.
The proof for $\K$ goes the same way.
\end{proof}


Convolution with mollifiers $\phi_n=\sum_{|k|\le1/r_n}z^k$
allows to embed hyperfunctions $\B(\T)$ into $\G_{H,r}(\T)$,
preserving the usual product of $\A_1(\T)$~\cite{DHPV2}.


\begin{proposition}[(Completeness)]{} Without assuming completeness of $\olra E$,
$\F_r(\ola E)$ is complete, and $\F_r(\ora E)$ is sequentially complete.
\end{proposition}

\begin{proof}
If ${(f^m)}_{m}$ is a Cauchy sequence in $\F_r(\ola E)$, there are 
increasing sequences ${(m_\mu)},{(n_\mu)}\in\N^\N$ such that
$
    \forall\mu\in\N,\, \forall\,k,\ell\ge m_\mu$, $
    \limsup\limits_{n\to\infty} p^\mu_\mu\p{f^k_n-f^\ell_n}^{r_n}<\frac1{2^\mu}
\,$
and more precisely
$
 \forall\,k,\ell\in[m_\mu,m_{\mu+1}] ~
    \forall n\ge n_\mu: p^\mu_\mu\p{f^k_n-f^\ell_n}^{r_n}<\frac1{2^\mu}
.\,$
\newcommand\mm{{\bar\mu}}
\newcommand\sm[1]{^{\smash{m_{#1}}}}
Let $ \mm(n)=\sup\set{\mu\mid n_\mu\le n}$, and consider the 
sequence
$    \bar f={(f\sm{\mm(n)}_n)}_n 
.\,$
Then we have 
$f^m\to\bar f$ in $\F_r(\ola E)$. 
Indeed, for $\veps$ and $p^{\mu_0}_\nu$ given,
take $\mu>\mu_0,\,\nu$ such that $\frac1{2^\mu}<\frac12\veps$.
As $p^\mu_\nu$ is increasing in both indices,
we have for 
$m\in[m_{\mu+s},m_{\mu+s+1}]$:
\\[.5ex]\centerline{$
    p^{\mu_0}_\nu(f^m_n-\bar f_n)^{r_n}
  \le p^\mu_\mu(f^m_n-f_n\sm{\mu+s+1})^{r_n}
    + \sum_{\mu'=\mu+s+1}^{\mm(n)-1}
        p^{\mu'}_{\mu'}(f\sm{\mu'}_n-f\sm{\mu'+1}_n)^{r_n}~.
$}\\[.5ex]
For $n>n_{\mu+s}$, one has 
$n\ge n_{\mm(n)}$, thus
$
   p^{\mu_0}_\nu(f^m_n-\bar f_n)^{r_n}
   < \sum
   _{\mu'=\mu+s}^{\mm(n)} \frac1{2^{\mu'}} < \frac2{2^\mu} < \veps
$.
\\
As $\F_r(\ola E) $ is a metrisable space, this implies completeness.
\\
For Cauchy nets in $\F_r(\ora E)$,
we use that for all $\mu$ there is $ \nu(\mu) $ such that
$ p_{\nu(\mu)}^\mu  \le  p_{\nu(\mu+1)}^{\mu+1}$ and
$
  p_{\nu(\mu)}^\mu \left( f_n^m - f_n^p \right) ^{r_n}
  <\veps_\mu,
$
where ${(\veps_\mu)}_\mu$ decreases to zero.
With this, we can prove the sequential completeness of $\F_r(\ora E)$
by the same arguments as above.
\end{proof}



\begin{remark}[(discreteness of induced topology)]{}
In \cite{DHPV2} we have shown that a net 
${(\delta_n)}_n\in\olra E^\N$
such that $\forall \psi\in\olra E,\ \smash{
\lim\nolimits_{n\to\infty}\int^{\vphantom j}
_{\R^s}\delta_n\cdot\psi=\psi(0)}$,
cannot be bounded in $\olra E$, under very weak assumptions.
%
From this we deduce that the topology of any algebra containing $\delta$ and
$\olra E$ \textbf{\emph{must}} induce the discrete topology on $\olra E$.
\end{remark}

\section{Families of scales and asymptotic algebras}\label{sect:scales}

\label{sequences}

We now generalize the growth conditions.
Consider a family $r=\p{r^m}_m$ of sequences $\p{r_n^m}_n$ decreasing to zero
as $n\to\infty$.
Suppose that either\\[1ex]
\centerline{(I) ~
$\forall m\in\N:r^m=O(r^{m+1})$,
 ~~ or ~~ (II) ~
$\forall m\in\N:r^{m+1}=O(r^m)$.
}
\begin{theorem}
Define \(
  \F_r(\olra E) = \bigcap\nolimits_{m\in\N} \F_{r^m}(\olra E) ,\
  \K_r(\olra E) = \bigcup\nolimits_{m\in\N} \K_{r^m}(\olra E) 
\)
in case (I),
resp.
\(
  \F_r(\olra E) = \bigcup\nolimits_{m\in\N} \F_{r^m}(\olra E) ,~
  \K_r(\olra E) = \bigcap\nolimits_{m\in\N} \K_{r^m}(\olra E)
\) in case (II).
Then again, $\G_r(\olra E)=\F_r(\olra E)\,/\,\K_r(\olra E)$ is an algebra.
\end{theorem}

\begin{proof}
Using $\K_{r^m}(\olra E)\cdot\F_{r^m}(\olra E)\subset\K_{r^m}(\olra E)$
and Theorem~\ref{equiv-scales}, it is easy to verify 
that in both cases, $\F_r(\olra E)$ is a subalgebra and $\K_r(\olra E)$ is an ideal thereof.
\end{proof}

\begin{example}
For 
$r^m_n=1$ if $n\le m$,  $0$ elsewhere,
we get Egorov-type algebras.
\end{example}


\begin{definition} Let \,$\mathbf{a} = \left( a_m: \N\to\R_+ \right) _{m\in\Z}$
be an asymptotic scale, \ie\
$\forall m\in\Z$, $a_{m+1}=o(a_m),~a_{-m}=1/a_m$ and $\exists M\in\Z:a_M=o(a_m^2)$.
The asymptotic algebra defined by $\,\mathbf a\,$ and a locally 
convex algebra $(E,\P)$ is the factor space 
\(\displaystyle
  \A_{(\mathbf{a})}(E,\P) =
  \frac {
  \set{ f\in E^\N\mid \exists m\in\Z~~\forall p\in\P:~ p(f)=O(a_m) }^{\!\vphantom l}
  }{
  \set{ f\in E^\N\mid \forall m\in\Z~~\forall p\in\P:~ p(f)=o(a_m) }
  }
.
\)
\end{definition}

\begin{example}
(i)  
$a_m(n)=n^{-m}$ leads to Colombeau type generalized algebras.\\
(ii) 
$a_m=1/\exp^m$ 
($m$-fold iterated  $\exp$ function) gives 
exponential algebras~\cite{ds}.%
\end{example}

\begin{theorem} For $r^m_n=
{|\log a_m(n)|}^{-1}$,
we have $\displaystyle\G_r(E,\P)=\A_{(\mathbf a)}(E,\P)$.
\end{theorem}

\begin{proof}
If $p(f_n)<C\,a_m(n)=C\,e^{1/r^m}$ (for $a_m>1$), then
$\limsup (p\circ f)^{r^m}<\infty$ and $f\in\F_r(E,\P)$.
Conversely, if
$\limsup (p\circ f)^{1/|\log a_{\bar m}|}<C$ then
$p\circ f \le {(a_m)}^{\log C}$, $(a_m,C>1)$.
Using the third property of scales,
$\exists M: p\circ f= o(a_M)$.
\\
Now consider $ \forall \bar m:p\circ f=o(a_{\bar m}) $.
Take $ m \in \N $. Then, for any $ q\in\N $,
there is $
\hat m$ such that $ a_{\hat m}=o({a_m}^q)$
and $p\circ f = o(a_{\hat m}) = o( {a_m}^q ) =
o((e^{-1/r_m})^q) = o((e^{-q})^{1/r^m}) $.
Therefore $\limsup(p\circ f)^{r^m}\le e^{-q}$,
and as $q$ was arbitrary, we have $\ultra f_{p,r^m}=0$,
thus $f\in\K_r(E,\P)$.
\\
Finally assume $ \forall \bar m: \limsup p(f_n)^{r^{\bar m}}=0 $,
\ie\ $\forall C>0:(p\circ f)^{1/|\log a_{\bar m}|}=o(C)$, thus
$p(f)=O(C^{|\log a_{\bar m}|}=O({a_{\bar m}}^{|\log C|})$.
Now for any $m$, let $\bar m=m+1$ and $C=1/e$.
Then $p(f)=O(a_{\bar m})=o(a_m)$, as required.
\end{proof}

\def\mfrac#1#2{\raise1ex\hbox{$#1$}\!\!\Big/\!\!\lower1ex\hbox{$#2$}}

A second kind of ``asymptotic'' algebras is of the form
\\[1ex]\centerline{$\displaystyle
  \A^{(\mathbf a)}(E,\P) = \frac{
  \set{ f\in E^\N \mid \forall\sigma<0 ~~ \forall p\in\P:~ p(f)=o(a_\sigma) }
}{
  \set{ f\in E^\N \mid \exists\,\sigma>0~~\forall p\in\P:~ p(f)=o(a_\sigma) }
}
$,}\\[1ex]
where $\mathbf{a}={(a_\sigma)}_{\sigma\in\R }$ is a scale
(\ie\ $\forall\sigma>\rho,~a_{\sigma}=o(a_{\rho})$, etc.),
indexed by a real number.
As the subalgebra is here given as intersection and the ideal as union of sets,
this case is not covered by the previous one.

\begin{proposition} For $r^m=\frac1{|\log a_{1/m}|}$, we have
$  \A^{(\mathbf a)}(E,\P) = \F'_r(\P)/\K'_r(\P)$, with
\(\,  \F'_r(\P)=\F'_r(E,\P) =
  \set{ f\in E^\N \mid \forall m\in\N~\forall p\in\P:
    \ultra{f}_{p,r^m}\le1 }
\,\) and\\
 \( \K'_r(\P)=\K'_r(E,\P)=
  \set{ f\in E^\N \mid \exists m\in\N~\forall p\in\P:
    \ultra f_{p,r^m}<1 }
\,.\)
\end{proposition}

\begin{example}
$ a_\sigma(n)=e^{-n\,\sigma}$ 
gives 
algebras with infra--ex\-po\-nen\-tial growth~\cite{infexp}, 
of particular interest for 
embeddings of periodic hyperfunctions.
\end{example}

\section{Functorial properties}\label{sect:functor}

A map $ \vp:\olra E\to\olra F$ obviously extends canonically
to $\G_r(\vp):\G_r(\olra E)\to\G_r(\olra F)$ if 
for all $f\in\F_r(\olra E)$ and $k\in\K_r(\olra E)$,
we have\\[1ex]
\centerline{$  \begin{array}[t]{rrrr}
  (F_1):& \vp(f)={(\vp(f_n))}_n\in\F_r(\olra F)~,
&
\text{ ~ and ~ }
  (F_2):& \vp(f+k)-\vp(f)\in\K_r(\olra F)~.
  \end{array}
$}
\begin{definition}
The \textbf{$r$--extension of a map }$\varphi:\olra E\to\olra F$
satisfying the above conditions $(F_1),(F_2)$, is defined as 
the map $
  \G_r(\vp) :\G_r(\olra E)\to\G_r(\olra F)$ such that $
  	[f]\mapsto\vp(f)+\K_r(\olra F)
$, 
where $f$ is any representative of $[f]=f+\K_{r}(\olra E)$.
\end{definition}

\def\Q{\mathcal Q}
\begin{example}\label{lin-cont-ext}
Linear mappings $\vp$ of locally convex vector spaces
$(E,\P)\to(F,\Q)$ are continuous iff
~$\forall q\in\Q~~\exists p
\in P~~\exists c>0~~\forall x\in E :~
  q(\vp(x)) \le c\,p
  (x) ~
$.\\
Then, $\forall f\in E^\N:\ultra{\vp(f)}_{q,r} \le \ultra f_{p,r}$,
whence $(F_1)$ and $(F_2)$, using linearity.
\end{example}

\def\uds#1{\underset{\displaystyle\llap(#1\rlap)}}
\def\rp#1{^+_{r^{#1}}}

We consider again a sequence of scales $(r^m)$ such that
$r^{m+1}\smash{\uds\ge\le} r^m$. 
Let us denote  $\F^+_{r^m} = \F_{r^m}(\R_+,|\cdot|)$ and
$\K^+_{r^m} = \K_{r^m}(\R_+,|\cdot|)$.

\begin{definition}
\def\iitem[#1]{\textbf{#1}}\def\uds#1{}%
In case $r^{m+1}\uds\ge\le r^m$,
an increasing map $g:\R_+\to\R_+$ is called
\iitem[$r$--moderate] iff
$    \forall\uds Mm\in\N~\exists\uds mM\in\N: g(\F_{r^m}^+)\subset\F_{r^M}^+$,
and
\iitem[$r$--compatible] iff it is continuous at 0 and
$    \forall \uds mM\in\N~\exists \uds Mm\in\N:h(\K^+_{r^m})\subset\K^+_{r^M}$.
\\In case $r^{m+1}\uds\le\ge r^m$, the definitions hold with
$\forall m~\exists M\leftrightarrow\forall M~\exists m$ exchanged.
\end{definition}

These notions allow to characterize maps that extend canonically to $\G_r$:

\begin{definition}\label{d:cont-r-temp}
A map $\vp:(E,\P)\to(F,\Q)$ is \textbf{continuously $r$--temperate} iff :
(a) there is an $
    r\text{--moderate function } g $ such that \\\mbox{}\hfill$
    \forall q\in\Q ~~ \exists p\in\P ~~
    \forall f\in E: q(\varphi(f)) \le g(p(f))
$,\\
and (b) 
there is an 
    $r$--moderate function  $g$ 
     and an 
    $r$--compatible function $h$\\\mbox{}\hfill {such that}
    $\forall q\in\Q ~~ \exists p\in\P ~~
  \forall f,
   k\in E:
    q(\varphi(f+k) - \varphi(f)) \le g(p(f)) \, h(p(k))
$.
\end{definition}

\begin{theorem}
Any continuously $r$--temperate map $\varphi$ extends canonically to
\(
    \G_r(\vp):\G_r(E,\P)\to\G_r(F,\Q) \,,
\)
and this canonical extension is continuous for the topologies induced
by $(\ultra\cdot_{p,r^m}{)}_{p\in\P,m\in\N}$ resp.
$(\ultra\cdot_{q,r^m}{)}_{q\in\Q,m\in\N}$.
\end{theorem}

\begin{proof}
Condition (a) of Definition~\ref{d:cont-r-temp} implies $(F_1)$,
and (b) 
gives $(F_2)$. We omit the straightforward calculations, a bit lengthy 
in view of the four cases to be treated~\cite{DHPV3}. 
Continuity of $\G_r(\vp)$ is obtained in the same way as $(F_2)$, 
replacing $p(f)\in\F\rp m$ by $\ultra f_{p,m}\le K$, and
$p(k)\in\K\rp m$ by $\ultra k_{p,m}\le\veps$.
\end{proof}


\section{Association in $r$--generalized algebras}\label{sect:association}

In several situations, e.g. when solving 
PDE, strong equality is impossible to obtain or not needed, and
approximation expressed by \emph{association\/} is sufficient.


\def\ur#1{\ultra{#1}_{|\cdot|,r}}

\begin{definition}\label{def:Ns}
Generalized numbers $[x],[y]\in\C_r$ are \textbf{associated} 
iff $x-y$ is a null sequence,
$[x]\approx [y]\iff x-y\in N=\set{x\in\C^\N\mid\lim x=0}$.
\label{D:N}
For $s\in\R$, they are $s$--associated, $[x]\stackrel s\approx[y]$, iff
$x-y\in N^{(s)}=\set{x\in\C^\N\mid x_n=o(e^{-s/r_n})}$.
\end{definition}

\begin{remark}%
(i) The definition is well-posed 
since  $\K\avi\subset N$.\\
(ii) We have $N^{(s)}=e_r^{-s}\,N$, where $e_r={(e^{\frac1{r_n}})}_n$ 
represents a positive unit of $\C_r$.
\\[.5ex]
(iii) \label{rem:unit-balls}%
All elements of the open unit ball are associated to zero, 
$ \ur x<1\ \Longrightarrow$\linebreak
$x\in N$, and $x\in N\impl \ur x\le1$, but
 $\frac1{r_n}\underset{n\to\infty}\to\infty$
also verifies $\ur{\smash{\frac 1r}}=1$.
\end{remark}


\def\Jass{\mathop{\smash{\underset\J\approx}}}
\begin{definition}
If $\J$ is an additive subset of $\F_r(\olra E)$ containing 
$\K_r(\olra E)$, 
two elements $F,G\in\G_r(\olra E)$ are \textbf{$\J$--associated}, $F \Jass G$, iff
$   
   F-G\in\J/\K_r(\olra E).
$
\end{definition}

\begin{proposition} If $\J$ is absolutely convex,
the relation $\Jass$ is stable under multiplication with elements of the closed unit ball.
\end{proposition}

Clearly, $\Jass$ is compatible with derivation iff $\J$ is stable under differentiation.



\def\Js{\J^{(s)}}
\begin{definition}
We call $F,G\in\G_r(E,\P)$ {\bf strongly associated}, $F\simeq G$,
iff $\forall p\in\P:d\pr{}(F,G)<1$.
For any $s\in\R$, \textbf{strong $s$--association} is defined as
\myeqn{
	{ F\overset{s}\simeq G }
	\iff \forall p\in\P:d\pr{}(F,G)<e^{-s}
	\iff \smash{F\underset{\Js}\approx G}
}
where $ \Js 
 = \big\{ f\in E^\N\mid\forall p\in\P:\ultra f_{p,r}<e^{-s} \big\} 
 =:B_{e^{-s}}^{(\P)}
$.
\end{definition}

\begin{remark} If $ F\overset{s}\simeq G$ for all $s>0$, then $F=G$,
because $\bigcap_{s>0}\Js=o_{\G_r(E,\P)}$.
\end{remark}


In order to define \textbf{weak association}, notice that a continuous bilinear form
$\lr<>{\cdot,\cdot}:\olra E\times\DD\to\C$
canonically extends to $\G_r(\olra E)\times\DD\to\C_r$ 
(cf. Example~\ref{lin-cont-ext}).
This allows to define,
for any convex subset $M$ of $\F_r(\C,|\cdot|)$ containing $0_{\C_r}$,
subspaces $\J$ of the form
$J_M=\big\{ f\in\olra E^\N\mid\forall\psi\in\DD: \lr<>{f,\psi}\in M \big\}$.

\def\Jss{\J_{(s)}}
\begin{definition}
For 
$M=N^{(s)}$ 
resp.\ $M=\smash{B_{e^{-s}}^{(|\cdot|)}}
$, association with respect to $J_M$ is called
\textbf{weak ($s,\DD'$)--} resp.\ \textbf{strong--weak ($s,\DD$)--association}
and is written $\smash{F\overset s{\underset\DD\approx}G}$
\rlap($\iff \forall\psi\in\DD:\lr<>{F-G,\psi}\stackrel s\approx0$),
resp. $ F\overset{s}{\underset\DD\simeq} G$
\rlap($\iff \forall\psi\in\DD:{|\lr<>{F-G,\psi}|}_r< e^{-s}$).
If $s=0$, it is omitted from notation.
\end{definition}

\begin{example}
In Colombeau's case, 
$[f]$, $[g]$ are weakly $(s,\D')$--associated iff
$n^s(f_n-g_n)\to0$ in $\D'(\Omega)$.
For ultradistributions and for periodic hyperfunctions,
with $\DD=\D^{(m)}$, 
$\DD=\D^{\{m\}}$ 
resp. $\DD=\A(\T)$, this is a new construction.
\end{example}
\begin{proposition}
Strong--weak $(s,\DD)$--association implies $(s,\DD')$--association, which conversely
implies strong--weak $(s',\DD)$--association only for all $s'<s$.
\end{proposition}
\begin{proof}
This follows from $\ultra x_r<e^{-s}\!\impl\! x\in N^{(s)}\!\impl\!\ultra x_r<e^{-s'}$ for $s'<s$,
while a counter-example for $s'=s$ can be built as in Remark~\ref{rem:unit-balls}.
\end{proof}

In a forthcoming paper, we explain in detail how these concepts of association
are 
useful in the context of regularity theory and microlocal analysis.

\end{document}